\documentclass[12pt]{article}

\usepackage{amssymb,amsthm,amsmath}
\usepackage{enumerate,bm,mathrsfs}
\usepackage{nicematrix}
\usepackage{tabularx}
\usepackage{tikz,graphicx}
\usepackage{tikz-network}
\usepackage{caption, subcaption}
\usepackage{tikz-3dplot,pgfplots}
\usetikzlibrary{math}
\usepackage{geometry}
\usepackage{hyperref}
\usepackage{cleveref}

\pgfplotsset{compat=1.10}
\usetikzlibrary{angles}

\newtheorem{theorem}{Theorem}[section]
\newtheorem{corollary}[theorem]{Corollary}
\newtheorem{lemma}[theorem]{Lemma}
\newtheorem{proposition}[theorem]{Proposition}

\newtheorem{assumption}[theorem]{Assumption}

\theoremstyle{definition}

\DeclareMathOperator*{\argmin}{arg\,min}
\providecommand{\abs}[1]{\left\lvert#1\right\rvert}
\providecommand{\norm}[1]{\left\rVert#1\right\rVert}
\providecommand{\set}[1]{ \left\{ #1  \right\}  }
\providecommand{\setb}[2]{ \left\{ #1 \ \middle| \  #2 \right\}  }
\providecommand{\innprod}[2]{\left\langle #1, #2 \right\rangle}
\providecommand{\parenth}[1]{\left( #1 \right) }
\providecommand{\restr}[1]{\left. #1 \right|}

\newcommand{\raisedchi}{\raisebox{\depth}{\(\chi\)}}

\providecommand{\dist}{\rho}
\providecommand{\ngd}{\Omega}
\providecommand{\vertices}{\mathcal{V}}
\providecommand{\uverts}{\mathcal{U}}
\providecommand{\kverts}{\mathcal{K}}
\providecommand{\vones}{\mathbf{1}}
\providecommand{\vzero}{\mathbf{0}}
\providecommand{\genlap}{L}
\providecommand{\lap}{\mathcal{L}}
\providecommand{\nlap}{\lap_{n}}  
\providecommand{\dlap}{\genlap_\ngd}
\providecommand{\rwlap}{\lap_{rw}}

\providecommand{\ungd}{U}
\providecommand{\ungdc}{\widetilde{U}}
\providecommand{\kngd}{K}
\providecommand{\kngdc}{\widetilde{K}}
\providecommand{\lag}{\raisedchi}
\providecommand{\ulag}{\lag_\uverts}
\providecommand{\klag}{\lag_\kverts}
\providecommand{\loclag}{\bar{\lag}}

\providecommand{\neighbors}[1]{\mathcal{N}_{#1}}

\providecommand{\interior}[1]{\mbox{int}\left(#1\right)}
\providecommand{\ebdry}[1]{\partial\left(#1\right)}
\providecommand{\vbdryi}[1]{\delta_i\left(#1\right)}
\providecommand{\vbdryo}[1]{\delta_o\left(#1\right)}

\providecommand{\wsp}{\ell_2(\vertices,\mu)}
\providecommand{\wspd}{\ell_2(\vertices,d)}
\providecommand{\wspsub}[1]{\ell_2(#1,\mu)}

\providecommand{\decaybase}{\eta}
\providecommand{\weightmatrix}{W}
\providecommand{\dmaxneighbors}{d_{\mathcal{M}}}
\providecommand{\innerrad}{R_\ngd}

\providecommand{\keywords}[1]{\noindent \textit{Keywords:} #1}

\providecommand{\msc}[1]{\noindent \textit{MSC:} #1}

\title{Locally supported, quasi-interpolatory bases for the approximation of functions on graphs
\thanks{This research was supported by grant DMS-1813091 from the National Science Foundation.}
\thanks{The authors would like to thank F.J. Narcowich and J.D. Ward at Texas A\&M University for the helpful conversations during the onset of this project.}
\thanks{Python code is available online at 
\href{http://arxiv.org/}{arXiv}.}
}

\author{
E. Fuselier\thanks{High Point University, High Point, NC.}
\and
J. P. Ward\thanks{North Carolina A\&T State University, Greensboro, NC.} 
}

\begin{document}

\maketitle
\\

\msc{05C90, 41A05, 41A15, 65D07.}

\begin{abstract}
Graph-based approximation methods are of growing interest in many areas, including transportation, biological and chemical networks, financial models, image processing, network flows, and more. In these applications, often a basis for the approximation space is not available analytically and must be computed. We propose perturbations of Lagrange bases on graphs, where the Lagrange functions come from a class of functions analogous to classical splines. The basis functions we consider have local support, with each basis function obtained by solving a small energy minimization problem related to a differential operator on the graph. We present $\ell_\infty$ error estimates between the local basis and the corresponding interpolatory Lagrange basis functions in cases where the underlying graph satisfies a mild assumption on the connections of vertices where the function is not known, and the theoretical bounds are examined further in numerical experiments.  Included in our analysis is a mixed-norm inequality for positive definite matrices that is tighter than the general estimate 
$\norm{A}_{\infty} 
\leq 
\sqrt{n} \norm{A}_{2}$.
\end{abstract}

\section{Introduction}

The analysis of functions on domains with a graph structure has become an important topic in many areas. In this setting one seldom has access to closed-form representations of functions forming the basis for a given approximation scheme, which means that if needed the corresponding basis functions must be computed. This can be a computationally wasteful task for moderately sized graphs, and an intractable problem for larger ones. This paper focuses on the issue of efficiently creating locally supported basis functions on graphs, where the bases are derived from a differential operator.
\\
\indent Some well-known tools for approximation on graphs are as follows. Positive definite basis functions and kernels for machine learning were considered in \cite{erb19,erb20}. Dictionaries of functions on graphs were treated in \cite{shuman20}. There are also a variety of signal processing techniques on graphs that may be of interest to the reader \cite{shuman13,shuman16, sun19}. More recently, splines and wavelets on circulant graphs were studied in \cite{kotzagiannidis19}. The bases that we consider here are connected to the early work on kernels and splines on graphs \cite{smola03}, as well as the analysis of variational splines on graphs \cite{pesenson08,pesenson09,pesenson10}. \\
\indent Our approach is rooted in the continuous domain theory, where thin plate splines or other kernels can be used to create a ``local Lagrange" basis that mimics an interpolatory Lagrange basis (see, for example \cite{fuselier13,hangelbroek18}). In this case, each local Lagrange basis function is essentially a kernel quasi-interpolant obtained by enforcing Lagrange interpolation conditions on a small neighborhood. However, the graph setting differs in a few ways, primarily due the availability of the kernel, which is the Greens function for some differential or pseudo-differential operator. On a continuous domain, the Greens functions are readily available with user-friendly formulas. On a graph, the Greens functions are costly to compute as the pseudo-inverse of a large matrix. To avoid the usual Greens function approach, we use energy minimization to compute the local Lagrange functions. However, let us point out that the direct analog of the continuous construction is of interest theoretically, and we plan to pursue that construction further in a later work.\\
\indent In addition to using only local information to create each local Lagrange basis function, we set the values outside of a local footprint to zero - a decision that is motivated by and justified in cases where the full Lagrange functions have exponential decay properties \cite{ward20}. The local construction and compact support has several obvious benefits. Being locally constructed, the local Lagrange basis can be computed in parallel. Compact support implies that one can store the basis functions using sparse data formats, and when a target function is to be approximated at new sites, only the local Lagrange functions close to the newly acquired data have to be updated. We show that, with reasonable assumptions on the graph, these proposed functions can be made arbitrarily close to the full Lagrange functions, and we also demonstrate the utility of the local Lagrange basis in numerical experiments.

Throughout the paper, we represent graphs by their adjacency and Laplacian matrices. We consider the standard Laplacian, which is a positive semi-definite matrix, as well as the normalized and random-walk variants that are also symmetric in an appropriately defined weighted space. The bounds that we present are based on the spectral theory for these operators, their pseudo-inverses, and their submatrices. Further, the estimates presented include constants that are based on local conditions of the graph and can be explicitly computed. Included in our analysis is a mixed-norm inequality for positive definite matrices (\Cref{cor:pos_def}) that is tighter than the general estimate $\norm{A}_{\infty} \leq \sqrt{n} \norm{A}_{2}$. 

The paper is organized as follows. Section \ref{sec:notation} includes preliminaries and the notation used for the graphs, graph operators, function spaces, and graph splines, as well as any assumptions made on the graphs. The main results of the paper are in Section \ref{sec:lagrange_relation}, where we establish that the local Lagrange functions approximate their corresponding Lagrange functions. Numerical demonstrations are provided in Section \ref{sec:experiments}, and a mixed-norm inequality for positive definite matrices, which we believe may be a new result and is of interest on its own, is included in the appendix.

\section{Notation and Preliminaries}
\label{sec:notation}

\tikzmath{\vertexsize = .2;}
\begin{figure}[ht!]
\centering
\begin{tikzpicture}
\Vertex[label=$v$,opacity=0,size=2*\vertexsize]{v}  
\Vertex[x=-1,y=0,opacity=0,size=\vertexsize]{A1}
\Vertex[x=-1.5,y=-1,opacity=0,shape=rectangle,size=\vertexsize]{B11}
\Vertex[x=-2,y=0,size=\vertexsize]{C111}
\Vertex[x=-2.5,y=-1,size=\vertexsize]{C112}
\Vertex[x=-1,y=-1.5,opacity=0,size=\vertexsize]{C113}

\Vertex[x=-.7,y=-.7,opacity=0,size=\vertexsize]{B12}
\Vertex[x=0,y=-1.2,opacity=0,size=\vertexsize]{C121}

\Vertex[x=1,y=0,opacity=0,size=\vertexsize]{A2}
\Vertex[x=1.1,y=1,opacity=0,shape=rectangle,size=\vertexsize]{B21}
\Vertex[x=.5,y=1.6,opacity=0,size=\vertexsize]{C211}
\Vertex[x=1.1,y=2,size=\vertexsize]{C212}
\Vertex[x=2,y=1.4,size=\vertexsize]{C213}

\Vertex[x=1.7,y=.5,opacity=0,shape=rectangle,size=\vertexsize]{B22}
\Vertex[x=2.2,y=-.2,size=\vertexsize]{C221}
\Vertex[x=2.6,y=.4,size=\vertexsize]{C222}

\Vertex[x=-.7,y=.8,opacity=0,size=\vertexsize]{A3}
\Vertex[x=-1.5,y=.6,opacity=0,shape=rectangle,size=\vertexsize]{B31}
\Vertex[x=-1.3,y=1.6,opacity=0,size=\vertexsize]{C311}
\Vertex[x=-2.6,y=.6,size=\vertexsize]{C312}

\Vertex[x=-.3,y=1.3,opacity=0,shape=rectangle,size=\vertexsize]{B32}
\Vertex[x=-.8,y=2.5,size=\vertexsize]{C321}
\Vertex[x=-.3,y=2,size=\vertexsize]{C322}
\Vertex[x=.2,y=2,size=\vertexsize]{C323}

\Vertex[x=1,y=-1,opacity=0,size=\vertexsize]{A4}
\Vertex[x=.5,y=-1.4,opacity=0,shape=rectangle,size=\vertexsize]{B41}
\Vertex[x=-.5,y=-2,size=\vertexsize]{C411}
\Vertex[x=.5,y=-2,size=\vertexsize]{C412}
\Vertex[x=1.5,y=-1.8,size=\vertexsize]{C413}

\Vertex[x=1.5,y=-1,opacity=0,shape=rectangle,size=\vertexsize]{B42}
\Vertex[x=3,y=-.7,size=\vertexsize]{C421}
\Vertex[x=3,y=-1.4,size=\vertexsize]{C422}

\Edge[](v)(A1)
\Edge[](v)(A2)
\Edge[](v)(A3)
\Edge[](v)(A4)
\Edge[](A1)(B11)
\Edge[](A1)(B12)
\Edge[](A2)(B21)
\Edge[](A2)(B22)
\Edge[](A3)(B31)
\Edge[](A3)(B32)
\Edge[](A4)(B41)
\Edge[](A4)(B42)
\Edge[](B11)(C111)
\Edge[](B11)(C112)
\Edge[](B11)(C113)
\Edge[](B12)(C121)
\Edge[](B21)(C211)
\Edge[](B21)(C212)
\Edge[](B21)(C213)
\Edge[](B22)(C221)
\Edge[](B22)(C222)
\Edge[](B31)(C311)
\Edge[](B31)(C312)
\Edge[](B32)(C321)
\Edge[](B32)(C322)
\Edge[](B32)(C323)
\Edge[](B41)(C411)
\Edge[](B41)(C412)
\Edge[](B41)(C413)
\Edge[](B42)(C421)
\Edge[](B42)(C422)
\end{tikzpicture}
\caption{Given a center $v$, we will choose some subset $\ngd\subset\vertices$ containing $v$ on which to compute the local Lagrange functions. In this example nodes within $\ngd$ are unshaded. Unshaded circles are in $\interior{\ngd}$, unshaded squares comprise $\vbdryi{\ngd}$, and the shaded circles make up $\vbdryo{\ngd}$. Interestingly, note that it is possible to have points in the interior of $\ngd$ that are farther away from the center than points on the inner boundary. 
}
\label{fig:graph_example}
\end{figure}

We let 
$\mathcal{G} 
= 
\set{
\vertices,
\mathbf{E}\subset \vertices \times \vertices, 
w,
\rho
}$ 
denote a finite, connected, weighted, simple graph with vertex set $\vertices$, edge set $\mathbf{E}$, symmetric weight function 
$w: \vertices\times\vertices \rightarrow \mathbb{R}_{\geq 0}$, 
and distance function 
$\rho: \vertices \times \vertices \rightarrow \mathbb{R}_{\geq 0}$. A value for $w$ is positive if and only if the corresponding vertices are connected by an edge. The distance function $\rho$ is positive on distinct pairs of vertices. In this paper, the weight of an edge is interpreted as the inverse of the length of the edge (specified by $\rho$). For non-adjacent vertices $x$ and $y$, the distance $\rho(x,y)$ is the length of the shortest path connecting $x$ to $y$. Lastly, we use the notation $x\sim y$ to denote that $x$ and $y$ are connected by an edge.

Given $\ngd\subset \vertices$, there are several associated sets of interest (see Figure \ref{fig:graph_example}). We define the interior of $\ngd$ to be
\begin{equation}
    \interior{\ngd} = \{x\in\ngd\,|\,\neighbors{x}\subset\Omega\}.
\end{equation}
where we use $\neighbors{x}$ to denote the union of $\{x\}$ and all vertices sharing an edge with $x$.\footnote{In the literature sometimes $\neighbors{x}$ or $\mathcal{N}(x)$ is used to denote the ``open" neighborhood of $x$, which only contains $x$ if $x$ is connected to itself. Here the symbol is defined so that $x$ is always a member of $\neighbors{x}$. This is sometimes called the ``closed" neighborhood of $x$.} We use $\neighbors{\Omega}$ to denote the union of $\neighbors{x}$ for all $x\in\Omega$. The edge boundary of $\ngd$, denoted by $\ebdry{\ngd}$, consists of all edges in $\mathbf{E}$ with exactly one vertex in $\ngd$. The collection of all vertices in $\ngd$ contained in an edge from $\ebdry{\ngd}$, denoted by $\vbdryi{\ngd}$, is the inner boundary of $\ngd$, and the outer boundary $\vbdryo{\ngd}$ is the set of all vertices in the complement $\ngd^c$ that are contained in an edge from $\ebdry{\ngd}$.

\subsection{The Graph Laplacian and Variants}
\label{subsec:laplacian_variants}

The graph Laplacian is given by $\lap = D - A$, where $A$ is the the adjacency matrix determined by the weights and $D$ is the diagonal matrix whose $j$th diagonal entry $d_j:=D_{j,j}$ is equal to the $j$th row sum of $A$. Two commonly used normalized versions of the Laplacian, the \emph{normalized} graph Laplacian $\nlap$ and the \emph{random walk} graph Laplacian $\rwlap$, are given by 
\begin{equation}
\nlap 
:= I - D^{-1/2}AD^{-1/2},
\quad 
\rwlap := I - D^{-1}A,
\end{equation}
and there are other normalizations in the literature 
\cite{chung1996combinatorial, johnson2007effectiveness}. 
Among other applications, the Laplacian and its variants are used to study the topology/geometry of the graph, to measure the energy or smoothness of functions of the vertices, and to construct learning and signal processing tools on graphs (learning kernels, filters, wavelets, etc.). 
Which Laplacian to use in a given situation has been an interesting topic of discussion in recent years \cite{boukrab2021random,johnson2007effectiveness,shuman13}.

We will use the symbol $\genlap$ to denote a \emph{general Laplacian}, which we define to be an operator $\genlap:\vertices\rightarrow\vertices$ that shares the same sparsity and sign structure of $\lap$, i.e. a value in the matrix $\genlap$ is nonzero if and only if the corresponding element in $\lap$ is nonzero, and corresponding entries share the same sign \cite{biyikoglu2007laplacian}\footnote{The definition in \cite{biyikoglu2007laplacian} is restricted to symmetric operators, but note that here we also allow asymmetric matrices.}. 
Moving forward we will refer to $\genlap$ as a ``general Laplacian" or simply ``the Laplacian," and we will call $\lap$ ``the standard Laplacian" to avoid confusion. Some results in this paper will be worked out for a general Laplacian when convenient to do so, but we will focus on the standard and random-walk Laplacians when more specificity is required.

\subsubsection{Weighted Spaces and the Connection Between $\rwlap$ and $\nlap$}
\label{subsec:weighted_spaces}

On the surface $\nlap$ and $\rwlap$ look very different. Note that $\nlap$ is symmetric while $\rwlap$ is in general not, and from the point of view of measuring smoothness, $\rwlap$ has the intuitive property that it annihilates constant vectors, i.e. 
$\rwlap \vones = \vzero$,
whereas $\nlap$ annihilates $D^{1/2}\vones$ \cite{chung97}. However, too much focus on these differences can hide fact that $\nlap$ and $\rwlap$ are merely different representations of the exact same operator.

To see this, we refer the curious reader to the presentation of the random walk Laplacian in \cite{grigor18}. While not symmetric in the usual sense, $\rwlap$ is symmetric with respect to the weighted inner product
\begin{equation*}
\innprod{f}{g}_{\ell_2(\vertices,d)}
:=
\sum_{v\in\vertices}f_v\,g_v\,d_v,
\end{equation*}
where $d_v$ is the degree of $v$ and $f_v$, $g_v$ denote the values of $f$ and $g$ at $v$ \cite[Lemma 2.4]{grigor18}. Note that the columns of $D^{-1/2}$ form a canonical orthonormal basis for this space. It can then be shown without much difficulty that $\nlap:\ell_2(\vertices)\rightarrow \ell_2(\vertices)$ is simply $\rwlap:\wspd\rightarrow\wspd$ expressed in the coordinate system with respect to this basis. Indeed, this relationship can be readily seen through the similarity identity
\begin{equation}\label{eq:similarityrwlapnlap}
\nlap = D^{-1/2}\rwlap D^{1/2},
\end{equation}
which also demonstrates the fact that these Laplacians share the same spectrum. Our preference here is to use the random walk form of this operator since it satisfies the familiar derivative property of annihilating constant functions, and so that we can make use of symmetry we will use $\ell_2(\vertices,d)$ as its domain.

So that we can address the random walk and standard Laplacians in a single framework, we will sometimes work on general weighted spaces. In what follows we let $\mu:\vertices\rightarrow \mathbb{R}_{>0}$ be a measure on the vertices, and let $\wsp$ denote the space of functions in the inner product defined by 
\begin{equation}\label{eq:weightedinnprod}
\innprod{f}{g}_{\wsp}
:=
\sum_{v\in\vertices}f_v\,g_v\,\mu(v).
\end{equation}

\subsubsection{The Dirichlet Laplacian}
\label{subsubsec:dirichelt_laplacian}

If $f$ is a function on $\vertices$ and $\ngd\subset\vertices$, we use $f_\ngd$ to denote the restriction of $f$ to $\ngd$. With that, the \emph{Dirichlet Laplacian} on $\ngd$ corresponding to $\genlap$, denoted $\dlap$, is defined as follows. Given $f:\ngd\rightarrow\mathbb{R}$, let $\widetilde{f}$ be the extension of $f$ to all of $\vertices$ by setting $\widetilde{f}$ to be zero on $\vertices\backslash\Omega$, then define
\begin{equation*}
\dlap f 
:= 
\restr{\left(\genlap \widetilde{f}\right)}_\ngd.
\end{equation*}
 Since in general the value of $\genlap f$ at $v\in\vertices$ only depends on the values of $f$ on $\neighbors{v}$, note that if $f$ is defined on all of $\vertices$ and $f$ vanishes on $\vbdryo{\ngd}$, then the Dirichlet Laplacian coincides with $\genlap$ in the sense that
\begin{equation}\label{eq:restrlapf}
  \dlap \left(\restr{f}_\ngd\right) 
  = 
  \restr{(\genlap f)}_\ngd,
\end{equation}
Note also that $\dlap$ is simply the submatrix determined by the rows and columns of $\genlap$ corresponding to $\ngd$. 

\subsection{Polyharmonic Splines on Graphs and Lagrange Functions}
\label{subsec:polyharmonic}

Consider a target function of interest 
$f:\vertices\rightarrow \mathbb{R}$ 
whose values are given on a subset 
$\kverts\subset\vertices$, 
which we refer to as the \emph{known} vertices. Our aim is to provide an efficient and accurate method to predict values of the target function on the \emph{unknown} vertices 
$\uverts := \vertices\backslash\kverts$. 
Given a known vertex $v\in\vertices$, a Lagrange function centered at $v$, which we will denote as $\raisedchi_v:\vertices\rightarrow \mathbb{R}$, satisfies\footnote{The reader should note that when the function in question is $\lag$, the vertex subscript refers to the center of the Lagrange function instead of its argument. This should not cause confusion.}
\begin{equation}
\lag_v(x) 
=
\begin{cases}
1, & x=v\\
0, & x\in\kverts\backslash\{v\}.
\end{cases}
\end{equation}
We may omit the subscript when specifying the center is unnecessary. A full collection of Lagrange functions gives rise to the interpolant 
$s_{f,\kverts}:\vertices\rightarrow \mathbb{R}$
\begin{equation}
 s_{f,\kverts}(x) 
 := 
 \sum_{v\in\kverts} f_v\lag_v(x).
\end{equation}

Here we study Lagrange functions from function spaces analogous to radial basis functions and classical polyharmonic splines on $\mathbb{R}^d$. Interpolation schemes based on these yield unique interpolants satisfying a minimal energy property, where the energy is measured via a ``native space" norm or semi-norm (see, for example \cite[Corollary 10.25]{wendland05}). As in \cite{ward20}, we will concern ourselves with Lagrange functions that minimize the energy measured by applying the Laplacian $L$. Specifically, $\lag_v$ will satisfy
\begin{equation}\label{fulllagviamin}
\lag_v = \argmin_f \norm{Lf}_{\wsp}
\quad
\text{subject to}
\quad
f_x
=
\begin{cases}
1 & x=v\\
0 & x\in \kverts\backslash\set{v}
\end{cases}.
\end{equation}

Unfortunately, closed-form representations for Lagrange 
bases
on graphs are unknown in most cases, and they can be prohibitively expensive to compute. However, in \cite{ward20} it was shown that in some cases Lagrange functions enjoy exponential decay as the distance between the evaluation point and the center grows. This suggests that a Lagrange function might be approximated well by using only information local to its center. We refer to a suitable approximation to $\lag_v$ obtained using only information in a small neighborhood around $v$ as an associated \emph{local Lagrange function}, and we denote a local Lagrange function by $\loclag_v$ (or simply $\loclag$ if specifying the center is unnecessary).

\subsubsection{Function and Laplacian indexing and properties}
\label{subsubsec:laplacian_index}

To help us study the properties of the local Lagrange function in what follows, we partition the functions and the Laplacian matrix into a block structure based whether a given vertex is unknown, known, or a member of the neighborhood $\ngd$ which will be used to define $\loclag$. The main sets of interest are: 
\[
\ungd
:=
\uverts \cap \ngd,
\quad\quad 
\ungdc:=\uverts \backslash \ngd,
\quad\quad 
\kngd:=\kverts \cap \ngd,
\quad \quad
\kngdc:=\kverts \backslash \ngd.
\]
Given a function $f:\vertices\rightarrow\mathbb{R}$, we will use set subscripts to denote the function restricted to that set, i.e. $f_\kverts$ will be the vector of known values of $f$, etc. Also, we express $L$ as

\vspace*{1cm}
\begin{equation}\label{eq:lapblocks}
L
=
\begin{pmatrix}
\genlap_{\uverts} & \genlap_{\kverts} \\
\end{pmatrix}
=
\begin{pNiceMatrix}[columns-width=0.75cm]
L_{\ungd,\ungd} 
& L_{\ungd,\ungdc} 
& L_{\ungd,\kngd}   
& L_{\ungd,\kngdc} 
\\
L_{\ungdc,\ungd} 
& L_{\ungdc,\ungdc} 
& L_{\ungdc,\kngd} 
& L_{\ungdc,\kngdc} 
\\
L_{\kngd,\ungd} 
& L_{\kngd,\ungdc} 
&  L_{\kngd,\kngd}  
& L_{\kngd,\kngdc} 
\\
L_{\kngdc,\ungd} 
& L_{\kngdc,\ungdc} 
& L_{\kngdc,\kngd} 
& L_{\kngdc,\kngdc} 
\\
\CodeAfter
\OverBrace[shorten,yshift=3pt]{1-1}{1-1}{\genlap_{\ungd} }
\OverBrace[shorten,yshift=3pt]{1-2}{1-2}{\genlap_{\ungdc} }
\OverBrace[shorten,yshift=3pt]{1-3}{1-3}{\genlap_{\kngd} }
\OverBrace[shorten,yshift=3pt]{1-4}{1-4}{\genlap_{\kngdc} }
\end{pNiceMatrix},
\end{equation}
where each $L_{A,B}$ is the block of $L$ whose elements are $L_{a,b}$ with $a\in A$, $b\in B$, and each $L_A$ is the sub-matrix of $L$ whose columns correspond to the set $A$\footnote{The sole exception to this notation is when the set in question is $\ngd$; the symbol $L_\Omega$ is reserved to denote the Dirichlet Laplacian.}. Using this notation, the Dirichlet Laplace operator on $\ngd$ can be written as
\begin{equation}\label{eq:dlapblocks}
\dlap
=
\begin{pmatrix}
L_{\ungd,\ungd}  & L_{\ungd,\kngd}  \\
L_{\kngd,\ungd}  & L_{\kngd,\kngd}  
\end{pmatrix}
.
\end{equation}

In what follows we may look at further partitions of some of the submatrices, taking advantage of the sparsity of $L$. For instance, since $\neighbors{\kngdc} \cap \ungd\subset \vbdryi{\ngd}$, for $x\in \kngdc$ and $y\in \ungd$ such that $y\notin\vbdryi{\ngd}$, we get $L_{x,y}=0$. Thus if we define $\ungd_i:=\vbdryi{\ngd}\cap\ungd$, we can write
\begin{equation}\label{eq:sparsity_inner_boundary}
\genlap_{\kngdc,\ungd} 
=
\begin{pmatrix}
\genlap_{\kngdc,\ungd_i} &  \genlap_{\kngdc,\ungd\backslash \ungd_i}
\end{pmatrix}
=
\begin{pmatrix}
\genlap_{\kngdc,\ungd_i} &  \vzero
\end{pmatrix}.
\end{equation}

\subsection{Assumptions on $\mathcal{G}$}
\label{sec:assumptions}
In \cite{ward20}, it was shown that under certain conditions the Lagrange functions decay exponentially with respect to the distance from the center. We conjecture that a similar decay holds in greater generality, and we plan to address this further in a future paper. Our main focus here is on the local Lagrange functions, so we separate the two issues by assuming decay of the Lagrange functions.
 
\begin{assumption}\label{as:exp_decay}
The Lagrange functions have exponential decay, i.e. there are constants $C,T>0$ and $0<\decaybase<1$ such that for any vertex $v\in\kverts$ and any $w\in \vertices$, we have
\begin{equation*}\label{eq:exp_decay}
\abs{\lag_v(w)}
\leq
C \decaybase^{T \rho(v,w)}.
\end{equation*}
\end{assumption}

The additional assumption that we rely upon is in terms of the connections at each $v\in\uverts$. Specifically, we require that each unknown vertex is connected to at least one known vertex, and we quantify this with a uniform bound on the percentage that the unknown weights contribute to the total degree at unknown vertices. 

\begin{assumption}\label{as:unknowndegree}
There is a constant $0\leq\alpha < 1$ such that for all $v\in\uverts$ we have
\[\sum_{x\in\uverts} A_{x,v} \leq \alpha\,d_v.\]
\end{assumption}

\section{Relating Lagrange and local Lagrange functions}
\label{sec:lagrange_relation}

In this section we discuss the construction of a local Lagrange function\footnote{The subscript of the associated center is omitted here to simplify the notation in what follows.} $\loclag$ and derive error estimates illustrating that for a sufficiently large neighborhood $\ngd$, $\loclag$ approximates $\lag$ well in the infinity norm. To do both we rely on the characterization of $\lag$ as the solution of the energy minimization problem in \eqref{fulllagviamin}. 

It will be convenient to consider $\genlap$ expanded in terms of unknown and known parts, as in \eqref{eq:lapblocks}. Since only $\lag_\uverts$ needs to be determined,  solving \eqref{fulllagviamin} is equivalent to solving 
\begin{equation}\label{eq:lst_sq_lagrange}
\lag_\uverts
=
\argmin_f \norm{ L_\uverts f + L_\kverts \lag_\kverts }_{\wsp}
\end{equation}
To define the local Lagrange function $\bar{\raisedchi}$, we solve essentially the same problem, but interpolation constraints are only imposed within a neighborhood $\ngd$ of the center vertex $v$ using the Dirichlet Laplace operator on $\ngd$. 
Considering the decomposition of $\genlap_\ngd$ in \eqref{eq:dlapblocks}, the local analogue of \eqref{eq:lst_sq_lagrange} we use to define the local Lagrange function within $\ngd$ is
\begin{equation}
\loclag_\ungd
=
\argmin_f 
\norm{
\begin{pmatrix}
L_{U,U} \\
L_{K,U}
\end{pmatrix} 
f 
+ 
\begin{pmatrix}
L_{U,K} \\
L_{K,K}
\end{pmatrix} 
\loclag_\kngd}_{\wspsub{\ngd}},
\end{equation}
where the $\wspsub{\ngd}$ norm is measured as in \eqref{eq:weightedinnprod}, but only summing over members of $\ngd$. Outside of the neighborhood $\Omega$, we simply define the local Lagrange function $\loclag$ to be zero. This convenient choice is justified by the fact that $\lag$ is known to decay exponentially away from its center. 

We will compare the Lagrange function $\lag$ to $\loclag$ via the normal equations of the least squares problems above. The normal equation for \eqref{eq:lst_sq_lagrange} is given by
\begin{equation}
L_{\uverts}^{\top} \weightmatrix L_{\uverts} \ulag
= 
-L_{\uverts}^{\top}\weightmatrix L_{\kverts} \klag,
\end{equation}
where $\weightmatrix$ is a diagonal matrix with $\weightmatrix_{v,v} = \mu(v)$. Using \eqref{eq:lapblocks} and expanding  in terms of entries inside and outside of $\ngd$ leads to
\begin{equation}\label{eq:lagrange_normal_equations}
\begin{pmatrix} 
\genlap^{\top}_{\ungd}\weightmatrix\genlap_{\ungd} & \genlap^{\top}_{\ungd}\weightmatrix\genlap_{\ungdc} \\ \genlap^{\top}_{\ungdc}\weightmatrix\genlap_{\ungd} &\genlap^{\top}_{\ungdc}\weightmatrix\genlap_{\ungdc} \end{pmatrix}
\begin{pmatrix}
\lag_{\ungd}  \\
\lag_{\ungdc} 
\end{pmatrix}
=
-\begin{pmatrix} 
\genlap^{\top}_{\ungd}\weightmatrix\genlap_{\kngd} & \genlap^{\top}_{\ungd}\weightmatrix\genlap_{\kngdc} \\ \genlap^{\top}_{\ungdc}\weightmatrix\genlap_{\kngd} &\genlap^{\top}_{\ungdc}\weightmatrix\genlap_{\kngdc} \end{pmatrix} 
\begin{pmatrix}
\lag_{\kngd}  \\ 
\vzero 
\end{pmatrix},
\end{equation}
where we have used the fact that $\lag_{\kngdc} = \vzero$. 
We are most interested in the part of this equation whose entries correspond to $\ngd$ (i.e. the top blocks), from which we get 
\begin{equation}\label{eq:define_lag}
\genlap^{\top}_{\ungd}
\weightmatrix\genlap_{\ungd}\lag_{\ungd} 
+ 
\genlap^{\top}_{\ungd}\weightmatrix
\genlap_{\ungdc}\lag_{\ungdc} 
=  
-
\genlap^{\top}_{\ungd}\weightmatrix
\genlap_{\kngd}\lag_{\kngd}.    
\end{equation}
We would like to compare these to the equations that define the local Lagrange function. 
Expanding into known and unknown blocks using \eqref{eq:dlapblocks}, and using the fact that 
$\loclag_\kngd = \lag_\kngd$ gives
\begin{equation}\label{eq:define_local_lag}
\left(
\genlap^{\top}_{\ungd,\ungd}\weightmatrix_{\ungd}\genlap_{\ungd,\ungd}
+
\genlap^{\top}_{\kngd,\ungd}\weightmatrix_{\kngd}\genlap_{\kngd,\ungd}
\right)
\loclag_\ungd 
= 
-\left(
\genlap^{\top}_{\ungd,\ungd}\weightmatrix_{\ungd}\genlap_{\ungd,\kngd}
+
\genlap^{\top}_{\kngd,\ungd}\weightmatrix_{\kngd}\genlap_{\kngd,\kngd}
\right)
\lag_{\kngd},
\end{equation}
where here we use the abbreviated notation $\weightmatrix_{\ungd}:=\weightmatrix_{\ungd,\ungd}$, etc. only for the diagonal weight matrix (this should not cause confusion). To compare \eqref{eq:define_lag} and \eqref{eq:define_local_lag}, we begin by focusing on the right side of \eqref{eq:define_lag}. By expanding in blocks using \eqref{eq:lapblocks}, we get
\begin{align}
\begin{split}
\genlap^{\top}_{\ungd}\weightmatrix
\genlap_{\kngd}
\lag_{\kngd} 
&=
\left(
\genlap^{\top}_{\ungd,\ungd}\weightmatrix_{\ungd}
\genlap_{\ungd,\kngd}
+
\genlap^{\top}_{\ungdc,\ungd}\weightmatrix_{\ungdc}
\genlap_{\ungdc,\kngd}
+
\genlap^{\top}_{\kngd,\ungd}\weightmatrix_{\kngd}
\genlap_{\kngd,\kngd}
+
\genlap^{\top}_{\kngdc,\ungd}\weightmatrix_{\kngdc}
\genlap_{\kngdc,\kngd}
\right)
\lag_\kngd
\\
&=
\left(
\genlap^{\top}_{\ungd,\ungd}\weightmatrix_{\ungd}\genlap_{\ungd,\kngd}
+
\genlap^{\top}_{\kngd,\ungd}\weightmatrix_{\kngd}\genlap_{\kngd,\kngd}
\right)
\lag_\kngd
+
\left(
\genlap^{\top}_{\ungdc,\ungd}\weightmatrix_{\ungdc}\genlap_{\ungdc,\kngd}
+
\genlap^{\top}_{\kngdc,\ungd}\weightmatrix_{\kngdc}\genlap_{\kngdc,\kngd}
\right)
\lag_\kngd.
\end{split}
\end{align}
Next, note that if the center vertex of the Lagrange function is within the interior of $\ngd$, it will not have any neighbors in $\ungdc$ or $\kngdc$. Thus under this mild assumption we get
\begin{equation}
\genlap_{\ungdc,\kngd} \lag_{\kngd}
= 
\genlap_{\kngdc,\kngd}\lag_\kngd 
= 
\vzero,
\end{equation}
which implies that the right sides of \eqref{eq:define_lag} and \eqref{eq:define_local_lag} are equal.

Similarly, we expand the first term on the left side of \eqref{eq:define_lag} to get
\begin{align}
\begin{split}
\genlap^{\top}_{\ungd}\weightmatrix\genlap_{\ungd}\lag_{\ungd} 
&=
\left(
\genlap^{\top}_{\ungd,\ungd}\weightmatrix_{\ungd}\genlap_{\ungd,\ungd}
+
\genlap^{\top}_{\ungdc,\ungd}\weightmatrix_{\ungdc}\genlap_{\ungdc,\ungd}
+
\genlap^{\top}_{\kngd,\ungd}\weightmatrix_{\kngd}\genlap_{\kngd,\ungd}
+
\genlap^{\top}_{\kngdc,\ungd}\weightmatrix_{\kngdc}\genlap_{\kngdc,\ungd}
\right)
\lag_\ungd\\
&=
\underbrace{
\left(\genlap^{\top}_{\ungd,\ungd}\weightmatrix_{\ungd}\genlap_{\ungd,\ungd}
+
\genlap^{\top}_{\kngd,\ungd}\weightmatrix_{\kngd}\genlap_{\kngd,\ungd}\right)
}_{B}\lag_\ungd
+
\left(
\genlap^{\top}_{\ungdc,\ungd}\weightmatrix_{\ungdc}\genlap_{\ungdc,\ungd}
+\genlap^{\top}_{\kngdc,\ungd}\weightmatrix_{\kngdc}\genlap_{\kngdc,\ungd}
\right)\lag_\ungd,
\end{split}
\end{align}
where we note that the matrix $B$ in the first term is the same as the matrix on the left side of the normal equation for the local Lagrange function in \eqref{eq:define_local_lag}. 
With this and the discussion above, we now see that the difference between the Lagrange and local Lagrange on $\ngd$ is
\begin{equation}\label{eq:loc_lag_err_weighted}
\loclag_{\ungd} - \lag_{\ungd} 
=  
B^{-1}\genlap^{\top}_{\ungd}\weightmatrix\genlap_{\ungdc}\lag_{\ungdc} 
+  
B^{-1}
\left(
\genlap^{\top}_{\ungdc,\ungd}\weightmatrix_{\ungdc}\genlap_{\ungdc,\ungd}
+
\genlap^{\top}_{\kngdc,\ungd}\weightmatrix_{\kngdc}\genlap_{\kngdc,\ungd}
\right)
\lag_{\ungd}
\end{equation}
Finally, we refine the error further to reveal its dependence on the boundary of $\Omega$ using \eqref{eq:sparsity_inner_boundary}, where the reader will recall that we used $U_i=\vbdryi{\ngd}\cap \ungd$. From this we deduce that 
\begin{equation*}
    \genlap_{\kngdc,\ungd} \lag_{\ungd} 
    = 
    \genlap_{\kngdc,\ungd_i} \lag_{\ungd_i},
\end{equation*}
Similarly, we have
\begin{equation*}
    \genlap_{\ungdc,\ungd} \lag_{\ungd} 
    = 
    \genlap_{\ungdc,\ungd_i} \lag_{\ungd_i},
\end{equation*}
and thus \eqref{eq:loc_lag_err_weighted} becomes 
\begin{equation}\label{eq:loc_lag_err_weighted_simplified}
\loclag_{\ungd}  - \lag_{\ungd}  
=  
B^{-1}\genlap^{\top}_{\ungd}\weightmatrix
\genlap_{\ungdc}
\lag_{\ungdc} 
+  
B^{-1}
\left(\genlap^{\top}_{\ungdc,\ungd}\weightmatrix_{\ungdc}\genlap_{\ungdc,\ungd_i}
+
\genlap^{\top}_{\kngdc,\ungd}\weightmatrix_{\kngdc}\genlap_{\kngdc,\ungd_i}\right)\lag_{\ungd_i}.
\end{equation}
In the next sections we will derive bounds on the error in the infinity norm in the random walk and standard Laplacian cases. 

\subsection{Lagrange Comparison for the Random Walk Laplacian}
\label{subsec:comparison_weighted_spaces}
If $\genlap$ is the random walk Laplacian, the matrix $W$ is simply the degree matrix $D$, and the error in \eqref{eq:loc_lag_err_weighted_simplified} becomes
\begin{equation}\label{eq:loc_lag_err_rw}
\loclag_{\ungd}  - \lag_{\ungd} 
=  
B^{-1}
\underbrace{
\genlap^{\top}_{\ungd}\,D\,\genlap_{\ungdc}
}_{(I)}\lag_{\ungdc} 
+  
B^{-1}
\underbrace{
\left(
\genlap^{\top}_{\ungdc,\ungd}\,D_{\ungdc}\,\genlap_{\ungdc,\ungd_i}
+
\genlap^{\top}_{\kngdc,\ungd}\,D_{\kngdc}\,\genlap_{\kngdc,\ungd_i}
\right)
}_{(II)}
\lag_{\ungd_i},
\end{equation}
where the matrix $B$ is given by
\begin{equation*}\label{eq:rwB}
    B 
    = 
    \genlap^{\top}_{\ungd,\ungd}\,D_{\ungd}\genlap_{\ungd,\ungd}
    +
    \genlap^{\top}_{\kngd,\ungd}\,D_{\kngd}\genlap_{\kngd,\ungd}.
\end{equation*}
Since we wish to bound the infinity norm of the error in \eqref{eq:loc_lag_err_rw}, and we desire bounds for the norms of some of the matrices therein. 

We begin with $(I)$. Note that the rows of this matrix are indexed by $\ungd$ and the columns are indexed by $\ungdc$. Given $v\in\ungd$ and $y\in\ungdc$, the $(v,y)^{th}$ element of $(I)$ is given by
\begin{equation*}
    \parenth{
    \genlap^{\top}_{\ungd}\,D\,\genlap_{\ungdc}
    }_{v,y} 
    = 
    \sum_{x\in\vertices}\genlap_{x,v}\,d_x\,\genlap_{x,y} 
\end{equation*}
Note that the sum index $x$ cannot simultaneously be equal to $v$ and $y$ because $\ungd$ and $\ungdc$ are disjoint. Thus we can split the sum into its $x=v$, $x=y$, and $x\sim v,y$ terms:
\begin{equation}\label{eq:Iexpandentry}
\parenth{
    \genlap^{\top}_{\ungd}\,D\,\genlap_{\ungdc}
    }_{v,y} 
    =
    \genlap_{v,v}\,d_v\,\genlap_{v,y} 
    + 
    \genlap_{y,v}\,d_y\,\genlap_{y,y} 
    + 
    \sum_{x\sim v,y}\genlap_{x,v}\,d_x\,\genlap_{x,y}.
\end{equation}
Recall that in this case the diagonal entries of $\genlap=\rwlap$ are 1, and the off-diagonal entries are
\begin{equation*}
\genlap_{x,y} 
= 
(\rwlap)_{x,y} 
= 
\frac{-1}{d_x}
A_{x,y},
\end{equation*}
where $A$ is the symmetric weighted adjacency matrix of the graph. Using this to expand further, we get
\begin{align*}
\begin{split}
    \parenth{
    \genlap^{\top}_{\ungd}\,D\,\genlap_{\ungdc}
    }_{v,y}   
    &= 
    -2 A_{v,y} 
    + 
    \sum_{x\sim v,y} 
    \frac{1}{d_x}A_{x,v}A_{x,y}. 
\end{split}
\end{align*}
Now using Assumption \ref{as:unknowndegree}, we get the following bound on the $v^{th}$ absolute row sum of (I)
\begin{align*}
\begin{split}
    \sum_{y\in\ungdc}\left|(\genlap^{\top}_{\ungd}\,D\,\genlap_{\ungdc})_{v,y}\right| 
    & \leq  
    2 \sum_{y\in\ungdc}A_{v,y} 
    + 
    \sum_{y\in\ungdc}\sum_{x\sim v,y}
    \frac{1}{d_x}A_{x,v}A_{x,y}. \\
    & =  
    2 \sum_{y\in\ungdc}A_{v,y} 
    + 
    \sum_{x\sim v} 
    A_{x,v}
    \frac{1}{d_x}\sum_{y\in\ungdc}
    A_{x,y}. \\
    & \leq  
    2 \alpha d_v 
    + 
    \sum_{x\sim v} A_{x,v}
    \leq 
    (2\alpha +1)d_v.
\end{split}
\end{align*}
It follows that 
\begin{equation}\label{eq:Ibd}
\|(I)\|_{\infty} \leq (2\alpha+1) d_M,
\end{equation}
where $d_M:=\max_{v\in\ungd} d_v$.

Now we focus on the absolute row sums of $(II)$, and note the similar structure of the two terms in $(II)$. Let $X$ and $Y$ be disjoint subsets of $\vertices$, and consider $\genlap^{\top}_{X,Y}\,D_{X}\,\genlap_{X,Y}$. Note that since $X$ and $Y$ are disjoint, we are only dealing with the off-diagonal entries of each Laplacian in the product. Since in this case $\genlap = \rwlap$, we get $L_{x,y} = -d_x^{-1}A_{x,y}$ on the off-diagonal. Thus for $v,y\in Y$, the $(v,y)$ entry of the product is 
\begin{equation*}
\parenth{\genlap^{\top}_{X,Y}\,D_{X}\,\genlap_{X,Y}}_{v,y} 
= 
\sum_{x\in X}\frac{1}{d_x}A_{x,v}\,A_{x,y}.
\end{equation*}
This gives that the $v^{th}$ row sum is
\begin{equation}\label{eq:IIrowsums}
\sum_{y\in Y} 
\sum_{x\in X}\frac{1}{d_x}A_{x,v}\,A_{x,y} 
=
\sum_{x\in X} A_{x,v}\frac{1}{d_x}\sum_{y\in Y}\,A_{x,y}. 
\end{equation}

Now we consider the expression $(II)$ from \eqref{eq:loc_lag_err_rw},
\begin{equation*}
    (II) 
    = 
    \genlap^{\top}_{\ungdc,\ungd}\,D_{\ungdc}\,\genlap_{\ungdc,\ungd_i}
    +
    \genlap^{\top}_{\kngdc,\ungd}\,D_{\kngdc}\,\genlap_{\kngdc,\ungd_i},
\end{equation*}
and note that replacing the set index $\ungd_i$ with the larger set $\ungd$ simply results in appended columns and therefore a possibly larger row sum. Thus given $v\in U$ we use \eqref{eq:IIrowsums} to get
\begin{align*}
    \begin{split}
        \sum_{y\in \ungd_i}(II)_{v,y} & \leq \sum_{x\in \ungdc} A_{x,v}\frac{1}{d_x}\sum_{y\in \ungd}\,A_{x,y} + \sum_{x\in \kngdc} A_{x,v}\frac{1}{d_x}\sum_{y\in \ungd}\,A_{x,y}. \\
        & \leq \sum_{x\in \ungdc} A_{x,v} + \sum_{x\in \kngdc} A_{x,v} = \sum_{x\in \vertices\backslash\ngd}A_{x,v} \leq d_v, \\
    \end{split}
\end{align*}
which yields
\begin{equation}\label{eq:IIbdrw}
\|(II)\|_{\infty} \leq d_M.
\end{equation}

Now we will treat $B^{-1}$, where $B$ is given in \eqref{eq:rwB}. Our estimates will require a bound on its minimum eigenvalue.
\begin{lemma}\label{lemma:lambda_min_weighted}
If $\genlap = \rwlap$ and the graph satisfies Assumption \ref{as:unknowndegree}, the matrix
\[B 
    = 
    \genlap^{\top}_{\ungd,\ungd}\,D_{\ungd}\genlap_{\ungd,\ungd}
    +
    \genlap^{\top}_{\kngd,\ungd}\,D_{\kngd}\genlap_{\kngd,\ungd}\]
satisfies    
\begin{equation*}
    \lambda_{min}(B)\geq (1 - \alpha)^2 d_m,
\end{equation*}
where $d_m:=\min_{x\in \ungd}d_x$. 
\end{lemma}

\begin{proof}
The proof will involve submatrices of $\genlap = \rwlap$ and the symmetric normalized Laplacian $\nlap$. To avoid too many subscripts, we will use the notation $\nlap(\ungd,\ungd)$ and $\rwlap(\ungd,\ungd)$ to denote the submatrices obtained by using the rows and columns  associated with $\ungd$. We begin by observing that $B$ is the sum of positive semi-definite matrices, so the result will follow once we establish that
\begin{equation*}
    \lambda_{min}
    \parenth{
    \rwlap^\top(\ungd,\ungd)
    \,
    D_{\ungd}
    \rwlap(\ungd,\ungd)
    }
    \geq 
    (1 - \alpha)^2 \min_{x\in\ungd} d_x.
\end{equation*}

Recall that $\nlap$ and $\rwlap$ are similar via \eqref{eq:similarityrwlapnlap}, and because 
$D^{-1/2}\nlap D^{1/2} = \rwlap$, 
then 
\[D_\ungd^{-1/2}\nlap(\ungd,\ungd) D_\ungd^{1/2} 
= 
\rwlap(\ungd,\ungd),\] 
so the submatrices corresponding to $\ungd$ are also similar and thus share the same spectrum. Using the above identity also yields
\begin{equation*}
\rwlap^{\top}(\ungd,\ungd)\,D_{\ungd}\rwlap(\ungd,\ungd) 
=  
D_\ungd^{1/2}
\nlap(\ungd,\ungd)^{\top}\nlap(\ungd,\ungd) 
D_\ungd^{1/2}.
\end{equation*}
Recall that $\nlap$ is symmetric, and combine the above observations along with the fact that 
$\lambda_{min}(FG)\geq \lambda_{min}(F)\lambda_{min}(G)$ if $F$ and $G$ are symmetric positive semi-definite to arrive at
\begin{align*}
\begin{split}
    \lambda_{min}
    \parenth{
    \rwlap^{\top}(\ungd,\ungd)\,D_{\ungd}
    \rwlap(\ungd,\ungd)
    }
    &\geq 
    \min_{x\in U}d_x \lambda_{min}(\nlap(U,U))^2 \\
    &= 
    d_m \lambda_{min}(\rwlap(U,U))^2. 
\end{split}
\end{align*}
To finish the proof, we us a Gersgorin estimate on the eigenvalues of $\rwlap(\ungd,\ungd)$. Recall that the diagonals of $\rwlap$ are $1$, and Assumption \ref{as:unknowndegree} implies that the sum of the off-diagonal elements of $\rwlap(\ungd,\ungd)$ are bounded by $\alpha$. Thus Gersgorin's Circle Theorem gives
\begin{equation*}
    \lambda_{min}(\rwlap(\ungd,\ungd) \geq (1 - \alpha),
\end{equation*}
which completes the proof.
\end{proof}

Before we state our main theorem, we need a definition. Given a neighborhood $\Omega$ with associated center $v\in\ngd$, we define the \emph{inner radius} of $\Omega$ with respect to $v$ to be
\[\innerrad:=\min_{x\in\vbdryi{\ngd}}\dist(v,x).\]
Thus for any $x$ outside of $\interior{\ngd}$, the bound in Assumption \ref{as:exp_decay} gives $|\lag_v(x)|\leq C \decaybase^{T\innerrad}$. 

\begin{theorem}\label{thm:comparsion_error_rw}
Suppose that the graph satisfies Assumptions \ref{as:exp_decay} and \ref{as:unknowndegree}, and consider the case where $\genlap = \rwlap$ in the $\ell_2(\vertices,d)$ inner product. If the local Lagrange function $\loclag$ centered at $v$ is constructed using the neighborhood $\Omega$, then the $\ell_\infty$ error between $\loclag$ and the associated full Lagrange function $\lag$ satisfies
\begin{equation*}
\norm{
\lag - \loclag 
}_{\ell_\infty}
\leq
C 
\parenth{\sqrt{n_U}+1}\frac{\parenth{\alpha + 1}}{(1-\alpha)^2} 
\frac{d_M}{d_m}
\decaybase^{T\innerrad},
\end{equation*}
where $n_U$ is the cardinality of $U\subset \ngd$, $\innerrad$ is the inner radius of $\ngd$ with respect to $v$, and the constants $C,T>0$, and $\decaybase<1$ are from Assumption \ref{as:exp_decay}.
\end{theorem}
\begin{proof}
First we focus on the error within $\ngd$. Recall the the error is identically zero on the known vertices. Using  \eqref{eq:loc_lag_err_rw}, \eqref{eq:Ibd}, and $\eqref{eq:IIbdrw}$, we get
\begin{align*}
\norm{
\lag - \loclag 
}_{\ell_\infty(\ngd)} 
& \leq 
\|B^{-1}\|_\infty\left(\parenth{2\alpha + 1} d_M 
\|\lag\|_{\ell_\infty(\ungdc)} 
+ 
d_M 
\|\lag\|_{\ell_\infty(\ungd_i)}\right),
\end{align*}
and since $\innerrad\leq \dist(v,x)$ for all $x$ outside of $\interior{\ngd}$, we use the exponential decay of the Lagrange function in \eqref{eq:exp_decay} to get
\begin{equation*}\label{eq:rwloclagerr1}
\norm{
\lag - \loclag 
}_{\ell_\infty(\ngd)} 
\leq 
C\|B^{-1}\|_\infty2\parenth{\alpha + 1} 
d_M \decaybase^{T\innerrad}.
\end{equation*}
Since $B$ is symmetric positive definite we can use Proposition \ref{prop:pos_def} (see the appendix) along with the lower bound on $\lambda_{min}(B)$ in Lemma \ref{lemma:lambda_min_weighted} to bound $\|B^{-1}\|_\infty$ via
\[
\|B^{-1}\|_\infty \leq \frac{\parenth{\sqrt{n_\ungd} + 1}}{2\lambda_{min}(B)} \leq \frac{\parenth{\sqrt{n_\ungd} + 1}}{2 (1 - \alpha)^2 d_m}. \]
This combined with \eqref{eq:rwloclagerr1} gives the result on $\ngd$. The error estimate on $\vertices\backslash\ngd$ follows from the exponential decay of $\lag$, the fact that the local Lagrange function is set to zero on $\vertices\backslash\ngd$, and
\[ 
1
\leq 
\parenth{\sqrt{n_U}+1}
\frac{\parenth{\alpha + 1}}{(1-\alpha)^2} 
\frac{d_M}{d_m}.
\]
\end{proof}
\noindent Observe that the ratio $d_M/d_m$ is a measure of (local) uniformity in the data, reminiscent of the \emph{mesh ratio} that sometimes appears in the continuous domain approximation literature. 

\subsection{Lagrange Comparison for the Standard Laplacian}
\label{subsec:unweighted}
If $\genlap$ is the standard Laplacian $\lap$ on $\ell_2$, the weight matrix $W$ is the identity, and the error in \eqref{eq:loc_lag_err_weighted_simplified} becomes
\begin{equation}\label{eq:loc_lag_err}
\loclag_{\ungd}  - \lag_{\ungd} 
=  
B^{-1}
\underbrace{
\genlap^{\top}_{\ungd}\genlap_{\ungdc}
}_{(I)}\lag_{\ungdc} 
+  
B^{-1}
\underbrace{
\left(
\genlap^{\top}_{\ungdc,\ungd}\,\genlap_{\ungdc,\ungd_i}
+
\genlap^{\top}_{\kngdc,\ungd}\,\,\genlap_{\kngdc,\ungd_i}
\right)
}_{(II)}
\lag_{\ungd_i},
\end{equation}
where the matrix $B$ is given by
\begin{equation}
    B 
    = 
    \genlap^{\top}_{\ungd,\ungd}\genlap_{\ungd,\ungd}
    +
    \genlap^{\top}_{\kngd,\ungd}\genlap_{\kngd,\ungd}.
\end{equation}
We can bound $\|(I)\|_{\infty}$ and $\|(II)\|_{\infty}$ using techniques similar to \Cref{subsec:comparison_weighted_spaces}. As before, the bounds are in terms of the maximum and minimum degrees over a subset of the unknowns, with the subtle difference that here the maximum is taken over a set that includes any unknowns on the outer boundary of $\ngd$.
\begin{lemma}\label{lemma:I_II_unweighted}
Assuming the graph satisfies Assumption \ref{as:unknowndegree}, if $\genlap$ is the standard Laplacian we get
\[    
\|(I)\|_{\infty} 
\leq 
(2\alpha + 1)\dmaxneighbors^2,
\quad\quad 
\|(II)\|_{\infty} 
\leq 
\dmaxneighbors^2,
\]
where 
$\dmaxneighbors
:=
\max_{x\in\neighbors{\Omega}\cap \uverts}
d_x$.
\end{lemma}
\begin{proof}
Let $v\in\ungd$ and $y\in\ungdc$. Using an expansion similar to \eqref{eq:Iexpandentry}, we have
\begin{equation*}
\parenth{
    \genlap^{\top}_{\ungd}\genlap_{\ungdc}
    }_{v,y} 
    =
    \genlap_{v,v}\genlap_{v,y} 
    + 
    \genlap_{y,v}\genlap_{y,y} 
    + 
    \sum_{x\sim v,y}\genlap_{x,v}\genlap_{x,y}.
\end{equation*}
Recall that in this case the diagonal entries of $\genlap=\lap$ give the vertex degree and the off-diagonal entries are those of $-A$, so 
\begin{align*}
\begin{split}
    \parenth{
    \genlap^{\top}_{\ungd}\genlap_{\ungdc}
    }_{v,y} 
    &= 
    -(d_v+d_y)A_{v,y} 
    + 
    \sum_{x\sim v,y} A_{x,v}A_{x,y}. 
\end{split}
\end{align*}
Now using Assumption \ref{as:unknowndegree}, we get the following bound on the $v^{th}$ absolute row sum of (I)
\begin{align*}
\begin{split}
    \sum_{y\in\ungdc}\left|(\genlap^{\top}_{\ungd}\genlap_{\ungdc})_{v,y}\right| 
    & \leq  
     d_v\sum_{y\in\ungdc}A_{v,y}+\sum_{y\in\ungdc}d_y A_{v,y} 
    + 
    \sum_{y\in\ungdc}\sum_{x\sim v,y}
    A_{x,v}A_{x,y} \\
    & \leq  
     \alpha d_v^2+\dmaxneighbors\sum_{y\in\ungdc}A_{v,y} 
    + 
    \sum_{x\sim v}A_{x,v}\sum_{y\in\ungdc} A_{x,y} \\
    & \leq  
     \alpha d_v^2+\alpha \dmaxneighbors d_v 
    + 
    \sum_{x\sim v}A_{x,v} d_x \\    
    & \leq  
     2\alpha \dmaxneighbors^2 
    + 
    \dmaxneighbors\sum_{x\sim v}A_{x,v}
     \leq  
     (2\alpha + 1)\dmaxneighbors^2 .
\end{split}
\end{align*}
To deal with (II), we follow steps similar to \eqref{eq:IIrowsums}, noting that the off-diagonal entries of $\genlap = \lap$ are those of $A$. Given $v\in\ungd$, the corresponding row sum satisfies
\begin{align*}
    \begin{split}
        \sum_{y\in \ungd_i}(II)_{v,y} & \leq \sum_{x\in \ungdc} A_{x,v}\sum_{y\in \ungd}\,A_{x,y} + \sum_{x\in \kngdc} A_{x,v}\sum_{y\in \ungd}\,A_{x,y} \\
        & \leq \sum_{x\in \ungdc} A_{x,v}d_x + \sum_{x\in \kngdc} A_{x,v}d_x
        \leq \dmaxneighbors\sum_{x\in \ungdc} A_{x,v} + \dmaxneighbors\sum_{x\in \kngdc} A_{x,v} \\
        &= \dmaxneighbors\sum_{x\in \vertices\backslash\ungd}A_{x,v} \leq \dmaxneighbors^2\\
    \end{split}
\end{align*}
which yields
\begin{equation*}
\|(II)\|_{\infty} \leq \dmaxneighbors^2.
\end{equation*}
\end{proof}

We will also need a bound on the minimum eigenvalue of $B$. 
\begin{lemma}\label{lemma:lambda_min_unweighted}
If $\genlap = \lap$ and the graph satisfies Assumption \ref{as:unknowndegree}, the matrix
\[B 
    = 
    \genlap^{\top}_{\ungd,\ungd}\genlap_{\ungd,\ungd}
    +
    \genlap^{\top}_{\kngd,\ungd}\genlap_{\kngd,\ungd}\]
satisfies    
\begin{equation*}
    \lambda_{min}(B)\geq (1 - \alpha)^2 d_m^2,
\end{equation*}
where $d_m:=\min_{x\in \ungd}d_x$. 
\end{lemma}

\begin{proof}
Let $v\in\ungd$. By the assumption, the 
absolute
sum of the off diagonals across the $v^{th}$ row of $\genlap_{\ungd,\ungd}$ is bounded by $\alpha d_v$. 
Thus by Gershorin's Theorem we get
\begin{equation*}
    \lambda_{min}(\genlap_{\ungd,\ungd})\geq (1 - \alpha)\min_{x\in \ungd}d_x.
\end{equation*}
Note that $B$ is the sum of symmetric semi-positive definite matrices, including $\genlap_{\ungd,\ungd}^T\genlap_{\ungd,\ungd}$. 
Since $\genlap=\lap$ is symmetric, so is $\genlap_{\ungd,\ungd}$. Thus we have 
\begin{equation*}
    \lambda_{min}(B)\geq\lambda_{min}(\genlap_{\ungd,\ungd}^{\top}\genlap_{\ungd,\ungd})= \lambda_{min}(\genlap_{\ungd,\ungd})^2\geq (1 - \alpha)^2\min_{x\in \ungd}d_x^2.
\end{equation*}
\end{proof}
Now one can repeat the arguments in the proof of Theorem \ref{thm:comparsion_error_rw} with the bounds derived in this section to get the following.
\begin{theorem}\label{thm:comparsion_error_unweighted}
Suppose that the graph satisfies Assumptions \ref{as:exp_decay} and \ref{as:unknowndegree}, and consider the case where $\genlap = \lap$ in the $\ell_2$ inner product. If the local Lagrange function $\loclag$ centered at $v$ is constructed using the neighborhood $\Omega$, then the $\ell_\infty$ error between $\loclag$ and the associated full Lagrange function $\lag$ satisfies
\begin{equation*}
\norm{
\lag - \loclag 
}_{\ell_\infty}
\leq
C 
\parenth{\sqrt{n_U}+1}\frac{\parenth{\alpha + 1}}{(1-\alpha)^2} 
\frac{\dmaxneighbors^2}{d_m^2}
\decaybase^{T\innerrad},
\end{equation*}
where $n_U$ is the cardinality of $U\subset \ngd$, $\innerrad$ is the inner radius of $\ngd$ with respect to $v$, and the constants $C,T>0$, and $\decaybase<1$ are from Assumption \ref{as:exp_decay}.
\end{theorem}

\section{Experiments}

\label{sec:experiments}

\begin{figure}
     \centering
\includegraphics[width=0.3\textwidth]{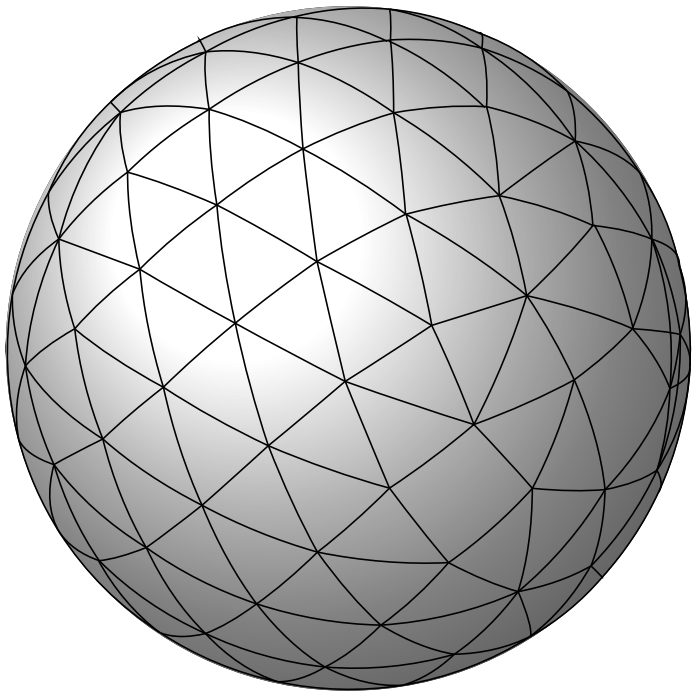}
\caption{ 
Graph generated from the Fibonacci lattice on the sphere with 125 points/vertices.
}
\label{fig:example_graph}
\end{figure}

In this section we illustrate the results presented above, with the aim of showing the promise of local Lagrange functions as computationally efficient bases on large data sets. First, we focus the convergence of the local Lagrange functions to the full Lagrange functions on graphs obtained by connecting points on the sphere $\mathbb{S}^2$. Next, we compare the computation time required to create the different bases at a fixed level of accuracy. Also, since we observed very similar results whether using the standard or random walk Laplacian, we only present the examples in the random walk case. The experiments were implemented in Python, and for large graphs, we use sparse data structures, i.e. SciPy sparse matrices. They were conducted on a desktop with 64GB of RAM and a Intel Xeon E5-1650 v3 CPU with frequency 3.50GHz and 6 cores.

In all of the experiments, we use data sites on the sphere $\mathbb{S}^2$ located at Fibonacci lattice points (\Cref{fig:example_graph}), where half of the vertices are unknown, and we choose the knowns as every other point in a standard Fibonacci lattice construction. This construction leads to graphs satisfying \Cref{as:unknowndegree}. In \Cref{subsec:connectivity}, we use a modification to this framework to explore cases where \Cref{as:unknowndegree} is not met. Since the points are nicely distributed, we use a subset to approximate the shortest geodesic distance between pairs of vertices, and connect pairs that are within a multiple (more precisely 2.5 in this case) of that minimum distance. Each edge weight is defined as the inverse of the geodesic distance between the vertices for the edge.

The neighborhoods $\ngd$ that define the footprint of the local Lagrange functions are determined by a maximum radius $R$: for each $v\in\kverts$ we define its neighborhood to be all vertices $x$ with $\rho(v,x)\leq R$, found using a single-source shortest path algorithm starting at the central vertex $v$ with the specified cutoff $R$. As before, we denote the number of unknowns in the neighborhood by $n_U$. These parameters specify a trade-off between computational cost and accuracy, and the experiments show that excellent results can be achieved with relatively small neighborhoods. After defining the neighborhoods, we compute the local Lagrange functions using a least-squares solver. 

Let us point out that the distances here (for non adjacent points) will be measured along the graph, not in the parameter space or on a manifold. However, for nicely distributed and dense data sets, these distances will be comparable.

\subsection{Convergence of bases and error bounds}

\begin{figure}
     \centering
\includegraphics[width=0.6\textwidth]{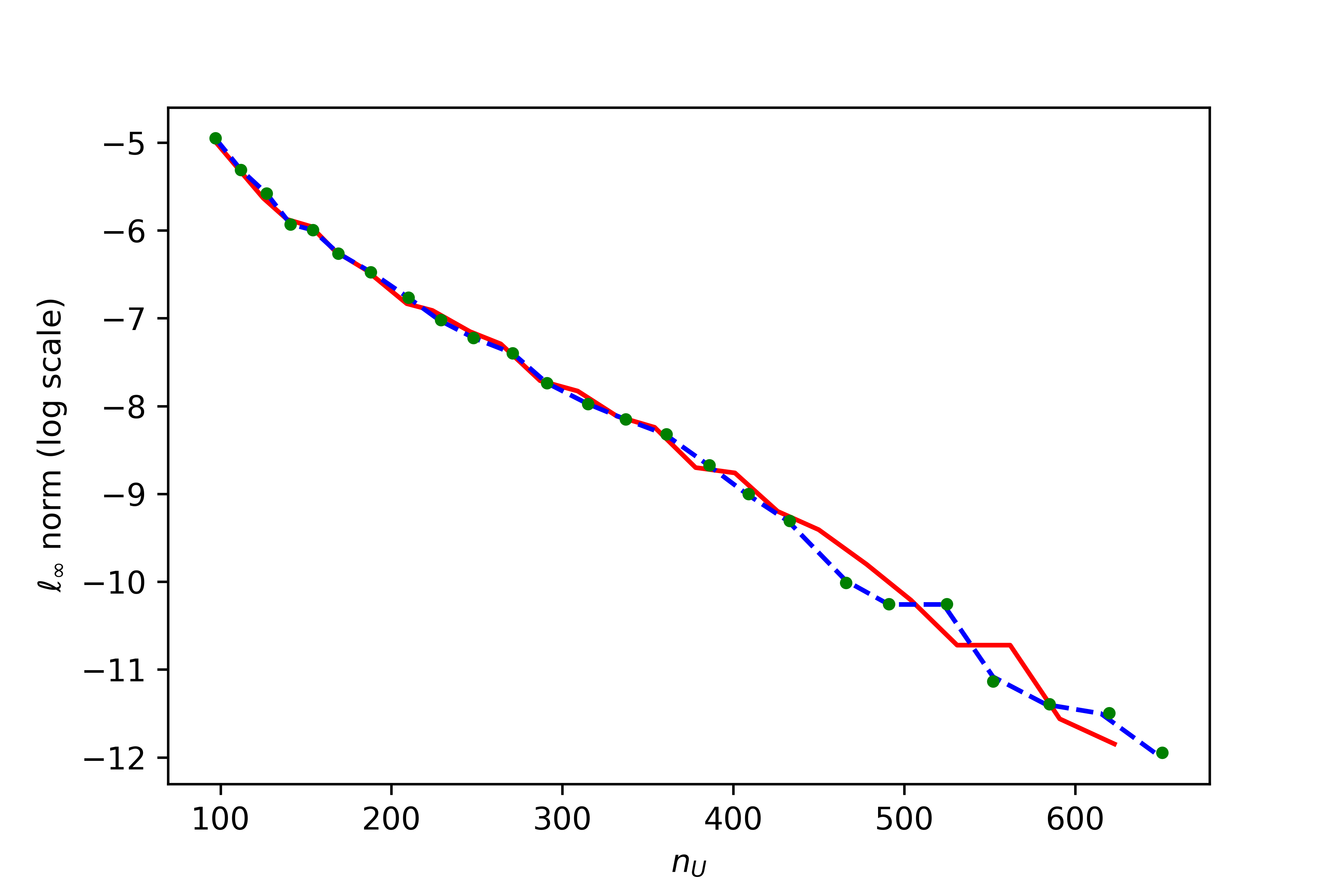}
\caption{ 
Comparison of a Lagrange and local Lagrange basis vector for 
$10^4$ (\textbf{solid red}),   
$10^5$ (\textbf{dashed blue}),
$10^6$ (\textbf{green dots})
points on the sphere.  The horizontal axis indicates the number of unknown points $n_U$ in the neighborhood used to define the local Lagrange function.
}
\label{fig:converge_error_sphere_single}
\end{figure}

Here through two examples we demonstrate convergence in the $\ell_\infty$ norm of the local Lagrange bases to the Lagrange basis with respect to the maximum support radius $R$ or the number of unknowns $n_U$ within the footprint. The first experiment illustrates convergence of a single local Lagrange basis function to the corresponding Lagrange basis function on large data sets of different sizes, from $10^4$ up to $10^6$ points on the sphere, (see \Cref{fig:converge_error_sphere_single}). Note that here roughly the same number of unknowns in the footprint achieves the same accuracy, independent of the total number of points on the graph. In this case, $n_U\sim 500$ consistently gives an error on the order of $10^{-10}$. In essence, this is a scaling property. 

\begin{figure}
     \centering
\includegraphics[width=0.6\textwidth]{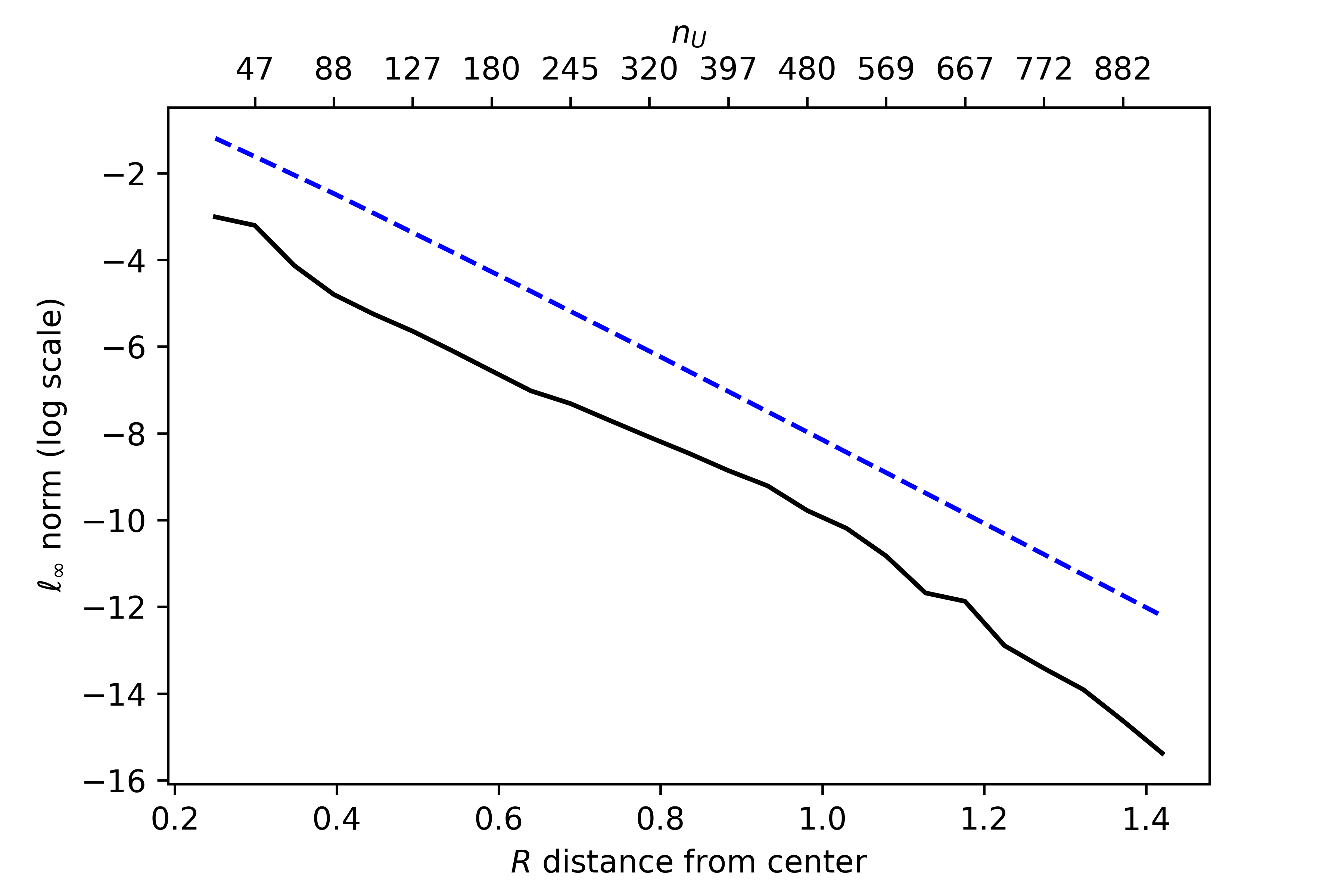}
\caption{ 
Illustration of the error bound using a graph on the sphere with $5,000$ vertices. The bottom horizontal axis denotes the radius $R$ of the footprint $\Omega$. The top horizontal axis indicates $n_U$, the number of unknown vertices inside $\Omega$.
\textbf{Dash Blue line:} 
error bound; 
\textbf{Black line:} 
Difference between Lagrange and Local Lagrange basis functions.
}
\label{fig:converge_error_sphere}
\end{figure}

In the second example, summarized in \Cref{fig:converge_error_sphere}, we focus on the bound in 
\Cref{thm:comparsion_error_rw}. We fixed the number of points to 5,000, computed the error between the local Lagrange and the Lagrange for varying footprints, and compared this with an estimate of the theoretical bounds\footnote{For simplicity, we use the maximum support radius $R$ instead of the inner radius $R_\ngd$ of each neighborhood to avoid the extra step of computing the interior of each footprint.}. The unknown constants $C$, $\eta$, and $T$ in the estimate were computed use a semilog best-fit line to the computed Lagrange decay, and the other constants in the bounds were computed explicitly. We note the consistent gap (a factor of about $100$) between the observed decay and the bound. This could illustrate that the estimate might be tightened, but the difference may also be due to the nice distribution of the Fibonacci lattice points, and there could be room for improvement in the bound by making additional assumptions on the data.

\subsection{Computation time}

\begin{figure}[ht]
    \centering
    \includegraphics[width=0.6\textwidth]{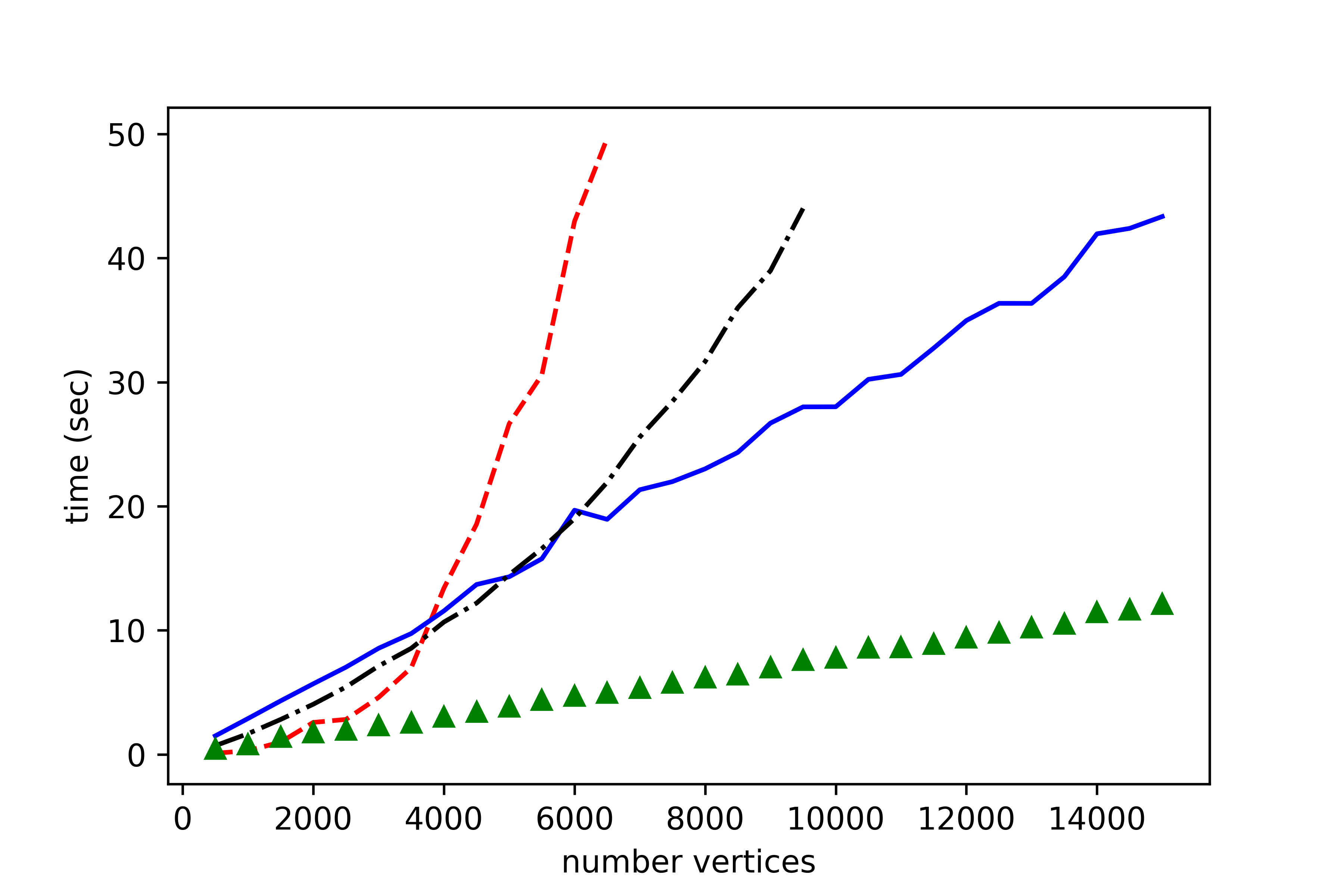}
    \caption{Computation time comparison between Lagrange and local Lagrange. \textbf{Red dashed line:} Lagrange basis, simultaneous computation in dense format;
    \textbf{Black dot-dashed line:} Lagrange basis, serial computation in sparse format; 
    \textbf{Solid blue line:} Local Lagrange basis, serial computation.
    \textbf{Green triangles:} Local Lagrange basis, parallel implementation.
    }
    \label{fig:compare_time}
\end{figure}

Using the same graphs on the sphere as in the previous section, we consider the time to compute the entire Lagrange and local Lagrange basis for various graph sizes, cf. \Cref{fig:compare_time}. For the full Lagrange basis, we consider two cases: simultaneous computation as the solution of a minimization problem and, taking advantage of the sparsity of the Laplacian, serial computation using a sparse solver. As the number of vertices grows so do the matrices in both cases, making the computation time increase faster than linearly, even when utilizing sparse solvers. For the local Lagrange basis, we choose neighborhoods giving $n_U \approx 150$, which consistently gave an error on the order of $10^{-5}$. 

First we computed the local Lagrange basis functions serially, and as expected we see a linear trend to the computation time, and we note that the time reported includes the construction of the neighborhoods. We also include the run times for a parallel version. A more robust parallel implementation of the local Lagrange basis would yield a drastic improvement in computation time. 

As a final point, note that given new data, updating the local Lagrange bases will also be comparatively fast, as this can be done locally. On the other hand, updating the Lagrange bases is equivalent to the original construction, due to global support of the basis functions.

\subsection{Connectivity of Unknowns}
\label{subsec:connectivity}

As a final experiment, we hint at future work by testing the necessity of Assumption \ref{as:unknowndegree}, which states that all unknown vertices be connected to at least one known vertex. We measured the error between a local Lagrange and corresponding Lagrange function on a graph constructed as before with $N=20,000$ nodes, but with varying percentages of known vertices: $20\%$, $15\%$, $10\%$, and $5\%$. The known vertices in each case were chosen in a way so that their distribution on the sphere was roughly uniform. 

\begin{figure}[ht]
    \centering
    \includegraphics[width=0.6\textwidth]{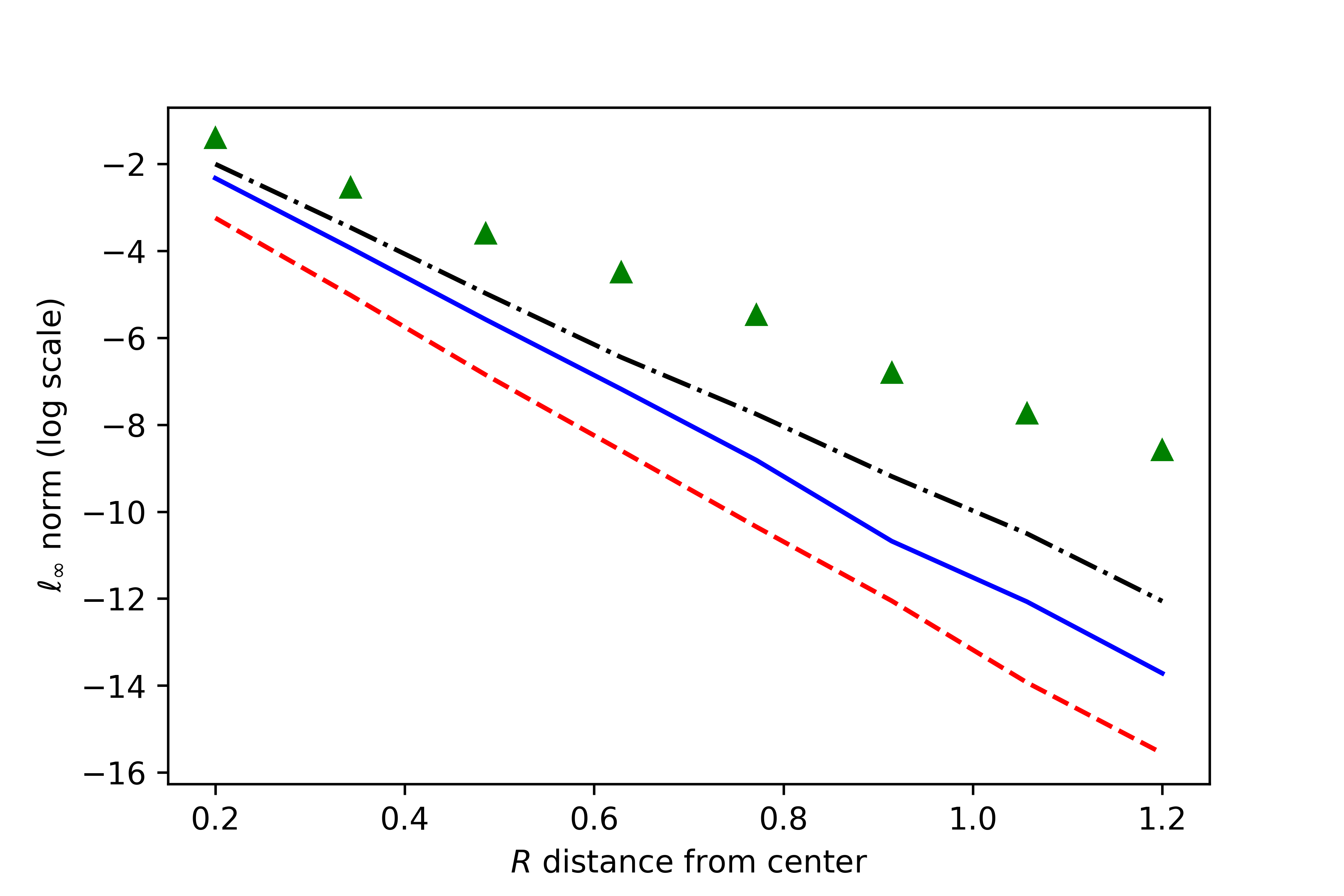}
    \caption{Comparison of a Lagrange and local Lagrange basis vector with  \textbf{Red dashed line:} $20\%$ known, $\alpha=0.94$, \textbf{Solid blue line:} $15\%$ known, $\alpha=0.95$, \textbf{Black dot-dashed line:} $10\%$ known, $\alpha=1$, and \textbf{Green triangles:} $5\%$ known, $\alpha=1$. The Horizontal axis represents the radius of the footprint used to compute the local Lagrange function.}
    \label{fig:low_perc}
\end{figure}

The $10\%$ and $5\%$ examples do not satisfy Assumption \ref{as:unknowndegree}. Indeed, in the $10\%$ case there were $221$ unknown vertices that were completely surrounded by unknowns, and in the $5\%$ case this figure rose to $4988$. Despite this, we see in \Cref{fig:low_perc} that the error between the local Lagrange and corresponding Lagrange basis function appears to decrease exponentially as the radius of the footprint $\Omega$ grows. We expect that the convergence of local Lagrange functions to the corresponding Lagrange function ultimately depends on the density of the knowns in the graph, and that Assumption \ref{as:unknowndegree} can be relaxed. 

\appendix

\section{Positive definite matrices}
\label{sec:pos_def}

It is possible that the following result is known. We were unable to find a reference, so we include a proof for the benefit of the reader.

\begin{proposition}\label{prop:pos_def}
If $A$ is an $n\times n$ positive definite matrix and $\lambda_{min}$ is the smallest eigenvalue of $A$, then
\begin{equation*}
\norm{A^{-1}}_{\infty}
\leq
\frac{\sqrt{n}+1}{2\lambda_{min}}.
\end{equation*}
\end{proposition}

\begin{proof}
The result will follow once we establish that
\begin{equation}\label{eq:pos_def_1}
\inf_{\norm{y}_{\ell_\infty}=1}
\norm{Ay}_{\ell_\infty}
\geq
\frac{2\lambda_{min}}{\sqrt{n}+1}.
\end{equation}
We begin by establishing that we may assume that the unit vectors $y$ in \eqref{eq:pos_def_1} are in the non-negative cone, i.e. $y_i\geq 0$ for all $i$. To see this, let $y$ be a vector that optimizes the left side of \eqref{eq:pos_def_1} and let $S$ be the diagonal matrix whose $i^{th}$ diagonal entry is 1 if $y_i\geq 0$ and $-1$ otherwise. Then $A_+ = SAS$ is symmetric positive definite, has the same eigenvalues of $A$, and it is not hard to show that a unit vector $y$ optimizes the left side of \eqref{eq:pos_def_1} if and only if the unit vector $Sy$, which is in the positive cone, optimizes the analogous quantity for $A_+$.

Next we proceed by writing $A = \lambda_{min} I +B $, where $B$ is positive semi-definite. For each $y$, let $By = \alpha_1y+\alpha_2\hat{y}$ where $\innprod{y}{\hat{y}}=0$. Since $B$ is positive semi-definite, $\alpha_1\geq 0$, and
\begin{equation}\label{eq:affinehyperplane_inf}
\norm{Ay}_{\ell_\infty}
=
\norm{\lambda_{min}y+\alpha_1y+\alpha_2\hat{y}}_{\ell_\infty}
\geq
\inf_{w\in \mathbb{R}^n,\innprod{y}{w}=0}
\norm{(\lambda_{min}+\alpha_1)y+w}_{\ell_\infty}.
\end{equation}
Note that the set
\begin{equation}
    P_y
    =
    \setb{(\lambda_{min}+\alpha_1)y+w}{\innprod{y}{w}=0}
\end{equation}
is an affine hyperplane. The crux of the next step can be seen in  \Cref{fig:inverseboundpf}, which shows that the smallest $\ell_\infty$ ball about the origin that intersects $P_y$ will contain a vector in the direction of $\vones$, where $\vones$ denotes the vector with all values being 1. Motivated by this, we solve for a point of the form
\begin{equation}
    \beta\vones 
    = 
    (\lambda_{min}+\alpha_1)y+w,\quad w\perp y,
\end{equation}
which yields
\begin{equation}
    \beta 
    = 
    \frac{(\lambda_{min}+\alpha_1)\norm{y}_{\ell_2}^2}{\norm{y}_{\ell_1}}.
\end{equation}
It is not hard to show (via contradiction) that indeed $\beta$ minimizes the $\ell_\infty$ norm for points in $P_y$. It now follows from \eqref{eq:affinehyperplane_inf} that
\begin{equation*}
\norm{Ay}_{\ell_\infty}
\geq
\norm{\beta \vones}_{\ell_\infty}
=
\beta.
\end{equation*}

\usetikzlibrary{calc}
\def\MarkRightAngle[size=#1](#2,#3,#4){
  \draw[black] ($(#3)!#1!(#2)$) -- ($($(#3)!#1!(#2)$)!#1!90:(#2)$) -- ($(#3)!#1!(#4)$)
}
\begin{figure}
  \begin{center}
\begin{tikzpicture}
  \coordinate (O) at (0,0);  
  \coordinate[label = above right:$(\lambda_{min}+\alpha_1)y$] (ay) at (1,5);
  \coordinate[label = below right:$y$] (y) at (3/5,3);
  \coordinate[label = above right:$\beta\vones$] (b) at (13/3,13/3);
  \MarkRightAngle[size=8pt](O,ay,b);
  \draw (0,0) node[below]{$\vzero$};
  \draw (0,3) node[left]{$\|x\|_{\ell_\infty} =1$};
  \draw (3,0) node[below]{$\|x\|_{\ell_\infty} =1$};
\draw[->] (0,0)--(4,0) node[right]{};
\draw[->] (0,0)--(0,4) node[above]{};
\draw[-] (3,0)--(3,3) node[above]{};
\draw[-] (0,3)--(3,3) node[above]{};
\draw[dashed,-] (13/3,0)--(b)--(0,13/3);
\draw[loosely dotted,thick] (O)--(ay);
\draw[loosely dotted,thick] (O)--(b);
\draw[<->,ultra thick] (-1,27/5)--(6,4);
\node at (y) [circle,fill,inner sep=1.5pt]{};
\node at (b) [circle,fill,inner sep=1.5pt]{};
\node at (ay) [circle,fill,inner sep=1.5pt]{};
\end{tikzpicture}
\caption{Here the bold line represents the affine hyperplane $P_y$ in the proof of \Cref{prop:pos_def}, motivating that a vector of the form $\beta\vones$ minimizes the infinity norm of vectors in $P_y$.}\label{fig:inverseboundpf}
\end{center}
\end{figure}

Next, since $\|y\|_{\ell_\infty} = 1$, for any $\ell_p$ norm we may write 
$\|y\|^p_{\ell_p} 
= 
1 + \|\bar{y}\|_{\ell_\infty}^p$ for some 
$\bar{y}\in \mathbb{R}^{n-1}$, allowing us to estimate further to get
\begin{equation*}
\norm{Ay}_{\ell_\infty}
\geq\beta =
\frac{(\lambda_{min}+\alpha_1)\norm{y}_{\ell_2}^2}{\norm{y}_{\ell_1}}
\geq
\lambda_{min}
\frac{\norm{y}_{\ell_2}^2}{\norm{y}_{\ell_1}}\\
\geq
\lambda_{min}
\inf_{\bar{y}\in\mathbb{R}^{n-1}}
\frac{1+\norm{\bar{y}}_{\ell_2}^2}{1+\norm{\bar{y}}_{\ell_1}}.
\end{equation*}
The greatest difference in $\|\bar{y}\|_{\ell_2}^2$ and $\|\bar{y}\|_{\ell_1}$ occurs in the direction $\bar{\vones}\in\mathbb{R}^{n-1}$, all entries being 1, which leads us to\footnote{This can be easily shown using Lagrange multipliers to solve: maximize $f = \sum_j x_j - \alpha^2$ subject to $\|x\|_{\ell_2}^2=\alpha^2$.}
\begin{equation*}
\inf_y \norm{Ay}_{\ell_\infty}
\geq
\lambda_{min}
\inf_{\bar{\beta}}
\frac{1+\norm{\bar{\beta}\bar{\vones}}_{\ell_2}^2}{1+\norm{\bar{\beta}\bar{\vones}}_{\ell_1}}
=
\lambda_{min}
\inf_{\bar{\beta}}
\frac{1+(n-1)\bar{\beta}^2}{1+(n-1)\bar{\beta}}.
\end{equation*}
Therefore, we want to minimize the function
\begin{equation*}
f(\bar{\beta})
=
\frac{1+(n-1)\bar{\beta}^2}{1+(n-1)\bar{\beta}}.
\end{equation*}
A basic calculus exercise shows that the minimum occurs at 
$\bar{\beta}=1/(\sqrt{n}+1)$ 
with the corresponding value for $f$ being $2/(\sqrt{n}+1)$. This establishes \eqref{eq:pos_def_1}.
\end{proof}
As an immediate consequence, we obtain the following mixed-norm inequality for positive semi-definite matrices, which are symmetric and have equal $\ell_1$ and $\ell_\infty$ norms. 
\begin{corollary}\label{cor:pos_def}
If $A$ is an $n\times n$ positive semi-definite matrix, then
\begin{equation*}
    \norm{A}_{1},\norm{A}_{\infty} \leq \frac{\sqrt{n}+1}{2}\norm{A}_{2}.
\end{equation*}
\end{corollary}

We can see that this is sharp by considering the matrix $B$ defined by 
\begin{equation*}
b_{i,j}
=
\begin{cases}
(\sqrt{n}+1)^2, & i=j=1\\
\sqrt{n}+1, & 
i\neq j \text{ and } (i=1 \text{ or } j=1)\\
1, & \text{otherwise}
\end{cases}
\end{equation*}
In this case,
\begin{align*}
\norm{B}_{\infty}
&=
\sqrt{n} (\sqrt{n}+1)^2\\
\norm{B}_{fro}^2
&=
(\sqrt{n}+1)^2
4n,
\end{align*}
\begin{align*}
\frac{\norm{B}_{\infty}}{\norm{B}_{2}}
\geq
\frac{\norm{B}_{\infty}}{\norm{B}_{fro}}
=
\frac{\sqrt{n}+1}{2}.
\end{align*}

\bibliographystyle{plain}
\bibliography{refs}

\end{document}